\documentclass{amsart}
\usepackage{epsfig}
\usepackage{mathrsfs}

\begin{document}

\newtheorem{thm}{Theorem}
\newtheorem{lem}[thm]{Lemma}
\newtheorem{cor}[thm]{Corollary}
\newtheorem{qn}{Question}

\theoremstyle{definition}
\newtheorem{defn}{Definition}

\theoremstyle{remark}
\newtheorem{rmk}{Remark}
\newtheorem{exa}{Example}

\def\square{\hfill${\vcenter{\vbox{\hrule height.4pt \hbox{\vrule width.4pt
height7pt \kern7pt \vrule width.4pt} \hrule height.4pt}}}$}

\def\R{\mathbb R}
\def\Z{\mathbb Z}
\def\CP{\mathbb {CP}}
\def\H{\mathbb H}
\def\E{\mathbb E}
\def\F{\mathscr F}
\def\D{\mathscr D}
\def\C{\mathbb C}
\def\til{\widetilde}
\def\N{\mathbb N}
\def\T{\mathscr T}
\def\G{\mathscr G}
\def\a{{\text{amb}}}

\newenvironment{pf}{{\it Proof:}\quad}{\square \vskip 12pt}

\title{Almost continuous extension for taut foliations}

\author{Danny Calegari}
\address{Department of Mathematics \\ Harvard \\ Cambridge, MA 02138}
\email{dannyc@math.harvard.edu}

\maketitle

\begin{abstract}
A taut foliation $\F$ of a hyperbolic $3$--manifold $M$ 
has the continuous extension property for
leaves in almost every direction; that is, for each leaf
$\lambda$ of $\til{\F}$ and almost every geodesic ray 
$\gamma$ in $\lambda$ the limit of $\gamma$ in $\til{M}$ is a well--defined
point in the ideal boundary of $\til{M} = \H^3$.
\end{abstract}

\vskip 12pt

\noindent{{\bf Acknowledgement:} I would like to thank the referee
for some useful comments.}

\vskip 12pt

Let $\F$ be a taut foliation of an atoroidal $3$--manifold $M$. Then a
theorem of Candel (\cite{aC93}) says that there is
a path metric on $M$ such that with their induced path metrics, leaves
of $\F$ are locally isometric to $\H^2$. In particular, it follows that
for any metric on $M$, the leaves of $\til{\F}$ with their induced path
metrics are uniformly quasi--isometric to $\H^2$, and therefore have a
well--defined {\em circle at infinity}. For a leaf $\lambda$ of $\til{\F}$
we denote this circle at infinity by $S^1_\infty(\lambda)$. Actually,
one only needs to know leaves of $\til{\F}$ are quasi--isometric to hyperbolic
planes to construct these circles at infinity, a fact which is much easier
to prove than Candel's theorem.

If $M$ is hyperbolic, there is an identification $\til{M} = \H^3$ and
there is a natural ideal boundary which we denote by $S^2_\infty$.

A basic problem in the theory of foliations of $3$--manifolds is to understand
the relationship between the intrinsic geometry of the leaves of the
foliation and the extrinsic (coarse) geometry of the ambient manifold (usually
in the universal cover). In particular, a question that has received a lot
of attention has been the following

\begin{qn}
Let $\F$ be a taut foliation of a hyperbolic $3$--manifold. Does the
inclusion $i:\lambda \to \til{M} = \H^3$ extend continuously to a map on
the ideal boundaries $i_\partial:S^1_\infty(\lambda) \to S^2_\infty(\H^3)$?
\end{qn}

This is a difficult question, and the (positive)
answer is only known in certain cases. In particular, the answer is known for
surface bundles over the circle (\cite{CT85}), for depth one and certain other
finite depth taut foliations (\cite{sF92},\cite{sF99}) and other special
cases. The problem is that leaves of taut foliations are far from being
quasi--isometrically embedded in the universal cover, so a path which is
quasi--geodesic in a leaf may potentially fail to limit to a definite point in
the ideal boundary of the ambient space. The subtlety of the
question is evidenced by the complicated structure of the image of such
$i_\partial$; limit sets of leaves are frequently ``exotic'' 
geometric sets such as dendrites, gaskets or sphere--filling curves.

It is easy to see, from the properness of leaves of $\til{\F}$, 
that if $i_\partial$ is defined for each point then it is continuous.
It turns out that if we only want to show that $i_\partial$ is defined
{\em almost everywhere}, then there is a surprisingly simple proof, which
works immediately for all taut foliations of hyperbolic $3$--manifolds.

\begin{thm}
With notation as above, for every leaf $\lambda$ of $\til{\F}$, for {\em almost
all} geodesic rays $r \subset \lambda$, the ray $r \subset \H^3$
converges to a definite point in $S^2_\infty$, and defines a measurable map
$i_\partial:S^1_\infty(\lambda) \to S^2_\infty(\H^3)$.
\end{thm}
\begin{pf}
Actually, the only property of a taut foliation we use is that $\lambda$
is quasi--isometric to $\H^2$ and the embedding $i:\lambda \to \H^3$ is
proper with bounded geometry and extends to a collar neighborhood of
$\lambda$. That is, there is an $\epsilon$ so that there is a quasi--isometric
embedding $I:\lambda \times [-\epsilon,\epsilon] \to \H^3$ such that
$I(p,0) = i(*)$.

That the geometry of the embedding is bounded follows from the compactness of
$M$. To see that there is a uniform collar 
neighborhood of $i(\lambda)$, observe
that there is a uniform $\delta$ so that the $\delta$--neighborhood of any
point in $M$ is contained in a ball foliated as a product. If $\lambda$ 
were to intersect such a lift of a product ball in two distinct disks in
$\til{M}$ we could find a transversal $\tau$ to $\til{\F}$ from $\lambda$ to
itself, and therefore by perturbation, a closed loop $\gamma \subset \til{M}$
transverse to $\til{\F}$. This contradicts the well--known fact that for a
taut foliation of a $3$--manifold,
transverse loops are homotopically essential.

Pick a basepoint $a \in \H^3$ through which the basepoint $p$ of $\lambda$
passes, and let $S$ be the visual sphere of $a$. There is an obvious
visual projection $\pi:\H^3 - a \to S$ which is basically just a version of
the Gauss map for hyperbolic space. Let $B$ be the ball of radius $1$ about $a$.

For a point $b \in \H^3$ with $\text{dist}_{\H^3}(a,b)=t$ and a vector
$v \in UT_b\H^3$ the norm $|d\pi(v)|$ is $O(e^{-t})$. The area of a
sphere of radius $t$ in $\H^3$ is $O(e^{2t})$. 
Therefore we can estimate 

$$\int_{\H^3 - B} \|d\pi(x)\|^\alpha d\text{vol}_{\H^3} 
\le \int_1^\infty \text{const.}\cdot e^{t(2-\alpha)} dt 
< \infty \text{ when } \alpha>2$$
where $\|\cdot\|$ denotes the operator norm on $d\pi:T\H^3 \to TS$.

Since $I$ is an embedding, we obviously have
$$\int_{I(\lambda \times [-\epsilon,\epsilon])\backslash B} 
\|d\pi(x)\|^\alpha d\text{vol}_{\H^3} < \infty$$
with the same assumption on $\alpha$, namely that $\alpha>2$.

The value of $\|d\pi(x)\|$ depends only on the distance from $x$ to $a$.
So, away from $B$, the fact that $I$ is a quasi--isometry implies that
this value varies only a bounded
amount over $q \times [-\epsilon,\epsilon]$ for each $q \in \lambda \backslash B$
{\em independently of $q$}, and so there is a bounded cost in replacing the
integral of this value over the interval with the value at $q \times 0$.

Moreover, we know $dI$ has uniformly bounded distortion over $\lambda$, by the
comments above, so that in particular
$d\text{vol}_\lambda dt \le \text{const.} I^*d\text{vol}_{\H^3}$. Thus

$$\int_{\lambda\backslash B} \|d\pi\circ i(x)\|^\alpha d\text{vol}_\lambda < \infty$$

In particular, using spherical co--ordinates on $\lambda$, we can conclude that
for almost every geodesic ray $\gamma \subset \lambda$ emanating from $p$,

$$\int_c^\infty \Biggl| \frac{d\pi(\gamma(t))} {dt} 
\Biggr|^\alpha e^t dt < \infty$$

so in particular, $e^{t/\alpha}\frac{d\pi(\gamma(t))}{dt}$ is in $L^\alpha$.
On the other hand, $e^{-t/\alpha}$ is certainly in $L^{\alpha/(\alpha-1)}$, so
by H\"older's inequality $\frac{d\pi(\gamma(t))}{dt}$ is in $L^1$.
That is, the position of
$\pi(\gamma(t))$ in $S$ moves only a bounded amount and therefore has a
well--defined limit. (Here we have chosen a hyperbolic metric on $\lambda$
according to Candel's theorem; such a metric is quasi--isometric to the
path metric on $\lambda$ inherited as a subspace of $\H^3$.)

It remains to show that $i_\partial$ is measurable. Recall that $\pi:\H^3 - a \to S$
is just radial projection onto the visual sphere at $a$. Let $\phi:S \to S^2_\infty$
be the visual identification. For any positive real number
$r$, let $S_r(p)$ denote the circle of radius $r$ in $\lambda$ centered
at $p$, and let $$f_r:S_r(p) \to S^1_\infty(\lambda)$$
be radial projection in $\lambda$.
Then $$\phi \pi f_r^{-1}:S^1_\infty(\lambda) \to S^2_\infty$$
is continuous for each $r$, and converges pointwise a.e. to $i_\partial$ as
$r \to \infty$. In particular, this limit is measurable.
\end{pf}

\begin{rmk}
Basically the point of the proof is the following: hyperbolic $3$--space
grows in volume like $O(e^{2t})$; a leaf of a taut foliation grows in
area like $O(e^{kt})$ where $-k$ is the ``coarse'' negative curvature
on a large scale of the leaf with respect to the induced subspace metric
(obviously $1 \le k \le 2$ and can be estimated from a
quasi--isometry constant of a uniformizing map $\lambda \to \H^2$). 
Since the embedding
of the leaf in $\H^3$ is a quasi--isometry into its $\epsilon$--neighborhood
for some $\epsilon$, the volume of a collar neighborhood of the leaf
can be efficiently measured by its area. Since the growth rate of both
$\lambda$ and $\H^3$ are exponential, it follows that ``most'' paths in
$\lambda$ make roughly comparable progress in $\H^3$ and in $\lambda$; that
is, ``most'' quasigeodesic rays in $\lambda$ limit to a definite 
point in $S^2_\infty$.
\end{rmk}

\begin{rmk}
A more subtle analysis of the properties of proper 
embedded minimal planes in $\H^3$ shows that we can actually
estimate $k \le \sqrt{2}$ (see \cite{wT97}).
\end{rmk}

\begin{rmk}
The proof applies essentially without modification to show that
leaves of the universal covers of essential laminations have measurable
extensions to $S^2_\infty$. The only technical issues are, firstly, that
essential laminations do not admit homotopically trivial tight transverse
loops (by \cite{GaOe}), and secondly that the leaves
of the universal cover of an essential lamination of a hyperbolic manifold,
with their induced path metrics,
are uniformly quasi--isometric to hyperbolic planes. This follows from
Candel's theorem.
\end{rmk}

\end{document}